\documentclass{elsart}
\usepackage{amsfonts,amssymb,amsmath}
\usepackage[notcite,notref]{showkeys}
\begin{document}
\begin{frontmatter}
\title{Gradings on simple algebras of finitary matrices}    
\author{Yuri Bahturin}\address{Department of Mathematics and Statistics\\Memorial
University of Newfoundland\\ St. John's, NL, A1C5S7, Canada\\
 and \\ Department of
Algebra, Faculty of Mathematics and Mechanics\\Moscow State University\\Moscow, 119992, Russia}\ead{yuri@math.mun.ca}\thanks{Work is partially supported by NSERC grant \#
227060-04 and URP grant, Memorial University of Newfoundland}
\author{Mikhail Zaicev}\address{Department of Algebra,
Faculty of Mathematics and Mechanics\\Moscow State University\\Moscow, 119992, 
Russia}\ead{zaicev@mech.math.msu.su}\thanks{Work is partially supported by RFBR, grant 06-01-00485a, and SSH-1983.2008.1}                   
\begin{abstract}
We describe gradings by finite abelian groups on the associative algebras of infinite matrices with finitely many nonzero entries, over an algebraically closed field of characteristic zero. 

\hfill\break 
 \hfill\break 
{\sl Mathematics Subject Index 2000:} 17B20, 17B40, 16W10, 16W50 

\end{abstract}

 \end{frontmatter}

\textwidth 140mm 
\textheight 210mm
\topmargin -0mm
\oddsidemargin 10mm
\evensidemargin 5mm
\newcommand{\pp}{\noindent {\em Proof. }}
\newtheorem{theorem}{Theorem}
\newtheorem{lemma}{Lemma}
\newtheorem{proposition}{Proposition}
\newtheorem{remark}{Remark}
\newtheorem{definition}{Definition}
\newtheorem{example}{Example}
\newcommand{\Aut}[1]{\mathrm{Aut}\,#1}
\newcommand{\bee}[1]{\begin{equation}\label{#1}}
\newcommand{\beq}[1]{\begin{eqnarray}\label{#1}}
\newcommand{\ben}{\begin{eqnarray*}}
\newcommand{\ene}{\end{equation}}
\newcommand{\eqe}{\end{eqnarray}}
\newcommand{\eqn}{\end{eqnarray*}}
\newcommand{\ld}{\ldots}
\renewcommand{\o}{\otimes}
\newcommand{\ra}{\rightarrow}
\newcommand{\diag}{\mathrm{diag}}
\newcommand{\cd}{\cdots}
\newcommand{\tr}[1]{\: ^t\!#1}
\newcommand{\rsl}{\mathrm{sl}}
\newcommand{\rgl}{\mathrm{gl}}
\newcommand{\wh}[1]{\widehat{#1}}
\newcommand{\wg}{\widehat{G}}
\newcommand{\wR}{\widetilde{R}}
\newcommand{\wC}{\widetilde{C}}
\newcommand{\wD}{\widetilde{D}}
\newcommand{\vp}{\varphi}
\newcommand{\bx}{\hfill$\Box$}
\newcommand{\cj}[2]{#1^{-1}\!#2#1}
\newcommand{\iv}[1]{#1^{-1}\!}
\newcommand{\ve}{\varepsilon}
\newcommand{\xba}{X_{\bar{a}}}
\newcommand{\xbb}{X_{\bar{b}}}
\newcommand{\xbab}{X_{\bar{a}\bar{b}}}
\newcommand{\xbe}{X_{\bar{e}}}
\newcommand{\xbg}{X_{\bar{g}}}
\newcommand{\ba}{\bar{a}}
\newcommand{\bz}{\bar{0}}
\newcommand{\bw}{\bar{1}}
\newcommand{\be}{\bar{e}}
\newcommand{\bg}{\bar{g}}
\newcommand{\bG}{\bar{G}}
\newcommand{\bh}{\bar{h}}
\newcommand{\bp}{\bar{p}}
\newcommand{\bq}{\bar{q}}
\newcommand{\bP}{\bar{P}}
\newcommand{\bQ}{\bar{Q}}
\newcommand{\bi}{\bar{i}}
\newcommand{\bj}{\bar{j}}
\newcommand{\bu}{\bar{u}}
\newcommand{\by}{\bar{y}}
\newcommand{\bv}{\bar{v}}
\newcommand{\ot}{\otimes}
\newcommand{\fd}{finite-dimensional\ }
\newcommand{\fdc}{finite-dimensional, }
\newcommand{\gh}{$\wg$-invariant}
\newcommand{\epf}{\hfill\qed}
\newcommand{\su}[1]{\mathrm{Supp}\,#1}
\newcommand{\ed}[1]{\mathrm{End}\,#1}
\newcommand{\aut}[1]{\mathrm{Aut}\,#1}
\newcommand{\Sp}[1]{\mathrm{Span}\,\{#1\}}
\newcommand{\twm}[4]{\left(\begin{array}{cc}#1&#2\\#3&#4\end{array}\right)}




\section{Introduction}\label{s0}

This paper is devoted to the extension of the results of \cite{surgrad,BZnc,BSZ} about the group gradings on finite-dimensional matrix algebras to the case of infinite-dimensional simple algebras of finitary linear transformations. After reminding the main results in the case of finite dimensions, we describe the $G$-graded embeddings of one finite dimensional graded matrix algebra into another (Theorem \ref{l992}), with $G$ a finite abelian group. Our next result says that if a simple algebra with minimal one-sided ideals is graded by $G$ as above then it can be presented as the direct limit of finite-dimensional $G$-graded matrix algebras (Theorem \ref{p9991}). This allows us to describe in Theorem \ref{t991} the gradings on the simple algebra of finitary matrices that is, the algebra of infinite matrices such that each matrix has only finitely many nonzero entries. Finally, in Theorem \ref{tu}, we give a necessary and sufficient condition for the equivalence of elementary gradings on the above algebra of infinite matrices.

\section{Some notation and simple facts}\label{s1}

Let $F$ be an arbitrary field, $R$ a not necessarily associative algebra over an $F$ and $G$ a group. We say
that  $R$ is a $G$-graded algebra, if there is a vector space sum
decomposition
\bee{e0001}
R=\bigoplus_{g\in G} R^{(g)},
\ene
such that
\bee{e0002}
R^{(g)}R^{(h)}\subset R^{(gh)}\mbox{ for all }g,h\in G.
\ene

Two $G$-gradings
\bee{e1001}
R=\bigoplus_{g\in G} R^{(g)},\mbox{ and }R=\bigoplus_{g\in G} (R^{\prime})^{(g)}
\ene
are called \emph{isomorphic} if there is an automorphism $\vp$ of $R$ such that $\vp(R^{(g)})=(R^{\prime})^{(g)}$, for all $g\in G$.

A subspace $V\subset R$ is called {\em graded} (or {\it homogeneous})
if $V=\oplus_{g\in G} (V\cap R^{(g)})$. An element $a\in R$ is called
{\it homogeneous of degree} $g$ if $a\in R^{(g)}$. We also write $\deg a=g$. The {\it support} of the
$G$-grading is a
subset
$$
\mathrm{Supp}\: R=\{g\in G|R^{(g)}\ne 0\}.
$$

\section{Reminder: Group gradings on matrix algebras}\label{stwotypes}

Below we briefly recall the results of \cite{surgrad,BZnc,BSZ}, where the full description of a finite group gradings on the full matrix algebra has been given.

\bigskip

A grading $R=\oplus_{g\in G}R^{(g)}$ on
the matrix algebra $R=M_n(F)$ is called {\it elementary} if there exists an $n$-tuple
$(g_1,\ldots,g_n)\in G^n$, such that the matrix units $E_{ij}, 1\le
i,j\le n$ are homogeneous and $ E_{ij}\in R^{(g)}\iff g=g_i^{-1}g_j. $ If $R$ is a matrix algebra with an elementary $G$-grading defined by a tuple $(g_1,\ldots,g_n)\in G^n$ and $B$ an algebra with a $G$-grading then the tensor product $R=A\otimes B$ will be given a grading if, given a homogeneous element $x$ of degree $h$, we set $E_{ij}\otimes x\in R^{(g)}$ provided that $g=\iv{g_i}hg_j$, for any $1\le i,j\le n$. This grading of the tensor product is called \emph{induced}.

A grading is called {\it fine} if $\dim R^{(g)}=1$ for any $g\in
\mathrm{Supp}\: R$. In this case $T=\mathrm{Supp}\: R$ is always a subgroup of $G$ \cite{BSZ}. In this case if $V$ is a natural $R$-module then $V$ is the space of a faithful irreducible representation of $T$ (see \cite{BZnc}). If we denote by $X_t$ the image of $t$ in $R$ corresponding to this representation then $X_t$ is a basis of $R_t$ and there is a $2$-cocycle $\alpha:G\times G\ra F^{\ast}$ such that $X_tX_s=\alpha(t,s)X_{ts}$, for any $t,s\in T$. This makes $R$ isomorphic to a twisted group algebra $F^{\alpha}G$.

The main result of \cite[Theorem 6]{surgrad} can be formulated as follow.

\begin{theorem}\label{tbz} Let $G$ be a group of order $d$, $F$ an algebraically
closed field  and $R=M_n(F)$. Then, as a $G$-graded algebra, $R$ is isomorphic to the tensor product with induced grading
$R\cong A\otimes B$
where $A=M_{k}(F)$ has an elementary $G$-grading, with support $S$, $B=M_{l}(F)$ has a fine grading, with support $T$,
and $S\cap T=\{e\}$.
\end{theorem}

A particular case of the fine gradings is a so-called
$\varepsilon$-grading where $\varepsilon$ is an $n^\mathrm{th}$ primitive root
of $1$. Let $G=\langle a\rangle_n\times\langle b\rangle_n$ be the direct product of two cyclic groups of order
$n$ and

\bee{e00}
X_a = \left(\begin{array}{cccc} \varepsilon^{n-1} & 0 & ... & 0
\\ 0 & \varepsilon^{n-2} & ... & 0
\\  \cdots& \cdots & \cdots& \cdots
\\ 0 & 0 & ... & 1
\end{array}\right)~,~~
X_b = \left(\begin{array}{cccc}
  0 & 1 &  ... & 0
\\ \cdots & \cdots   &    \cdots &  \cdots
\\ 0 & 0 &  ... & 1
\\ 1 & 0 &  ... & 0
\end{array}\right)~.
\ene
Then
\begin{equation}\label{a1}
X_aX_bX_a^{-1}=\varepsilon X_b~,~~X_a^n=X_b^n=I
\end{equation}
and all $X_a^iX_b^j, 1\le i,j\le n$, are linearly independent.
Clearly, the elements $X_a^iX_b^j,\, i,j=1,\ldots, n$, form a basis
of $R$ and all the products of these basis elements are uniquely
defined by (\ref{a1}).

       Now for any $g\in G, g=a^ib^j$, we set $X_g=X_a^iX_b^j$ and denote by $R^{(g)}$
a one-dimensional subspace
\bee{aa1}
R^{(g)}=\langle X_a^iX_b^j\rangle.
\ene
Then from (\ref{a1}) it follows that $R=\oplus_{g\in G}R^{(g)}$ is a
$G$-grading on $M_n(F)$ which is called an $\varepsilon$-grading.

Now let $R=M_n(F)$ be the full matrix algebra over $F$ graded by an abelian
group $G$. The following result has been proved in \cite[Section 4, Theorems 5, 6]{BSZ} and \cite[Subsection 2.2, Theorem 6, Subsection 2.3, Theorem 8]{surgrad}.

\begin{theorem}\label{tbsz} Let $F$ be an algebraically closed field of characteristic zero. Then as a $G$-graded algebra $R$ is isomorphic to the
tensor product
$$
R_0\otimes R_1\otimes\cd\otimes R_k
$$
where $R_0=M_{n_0}(F)$ has an elementary $G$-grading, $\mathrm{Supp}\:
R_0=S$ is a
finite subset of $\,$ $G$, $R_i=M_{n_i}(F)$ has the $\varepsilon_i$ grading,
$\varepsilon_i$ being a primitive $n_i^\mathrm{th}$ root of $1$,
$\mathrm{Supp}\; R_i=H_i\cong \mathbb{Z}_{n_i}\times
\mathbb{Z}_{n_i}, i=1,\ld,k$. Also $H=H_1\cd H_k\cong H_1\times\cd\times
H_k$ and $S\cap H=\{e\}$ in $G$.
\end{theorem}

\section{Embeddings of graded matrix algebras}
To describe the gradings on the algebra of finitary matrices we will need to consider the embedding of one $G$-graded \fd simple algebra into another.  We recall that if $R\cong M_n(F)$  is $G$-graded then, as a graded algebra, $R$ is isomorphic to a tensor product $C\ot D$ where $C=M_p(F)$, $D=M_q(F)$, $n=pq$, $M_p(F)$ is a matrix algebra with elementary $G$-grading, $M_q(F)$ a matrix algebra with a fine $T$-grading where $T=\su{D}$ is a subgroup in $G$ such that $T\cap\su{C}=\{ 1\}$ where $1$ is the identity element of $G$. Thus we may think that $R=CD\cong C\ot D$. Let us notice that the subalgebra $D$ is not defined uniquely and once $D$ has been fixed, $C$ is uniquely defined as the centralizer of $D$ in $R$.

\begin{theorem}\label{l992}
Let $R_1\cong M_k(F)$ and $R_2\cong M_n(F)$ be two $G$-graded matrix algebras with identity elements $e_1$ and $e_2$, respectively, $R_1=C_1D_1$, $R_2=C_2D_2$ their decompositions in which $C_1, C_2$ have elementary grading while $D_1, D_2$ have fine grading. Let also $D_1\cong D_2$ as graded algebras and $\vp:R_1\ra R_2$ be an injective homomorphism of graded algebras. Then there exists a decomposition $R_2=\wC_2\wD_2\cong\wC_2\ot\wD_2$ such that $\wC_2$ is a matrix algebra with elementary grading, $\wD_2$ is a matrix algebra with fine grading, $\wD_2\cong D_2$ as graded algebras and $\vp(C_1)\subset\wC_2$. If $\vp(e_1)R_2\vp(e_1)=\vp(R_1)$ then also $\vp(e_1)\wC_2\vp(e_1)=\vp(C_1)$. Besides, there is an isomorphism $\psi: D_1\ra\wD_2$ such that $\vp(a)\vp(d)=\vp(a)\psi(d)$, for any $a\in R_1$, $d\in D_1$.
\end{theorem}

\pp Since $D_1\cong D_2$, in particular, $\su{D_1}\cong\su{D_2}$. We denote $T=\su{D_1}$. Then for any $t\in T$ there exist invertible matrices $X_t\in R_1$ and $X_t^{\prime}\in R_2$ such that
$$
D_1=\Sp{X_t\:|\:t\in T},\;D_2=\Sp{X_t^{\prime}\:|\:t\in T}.
$$
In addition, for any $t,s\in T$ there is $\alpha(t,s)\in F$ such that
\bee{e995}
X_tX_s=\alpha(t,s)X_{ts},\; X_t^{\prime}X_s^{\prime}=\alpha(t,s)X_{ts}^{\prime},
\ene
following because $D_1$ and $D_2$ are isomorphic. One may also assume that $X_1$ and $X_1^{\prime}$ are the identity elements of $R_1$ and $R_2$, respectively.

Since $T\cap\su{C_2}=\{ 1\}$ it follows that $\vp(X_t)=A_tX_t^{\prime}$ for some matrix $A_t\in C_2$, $\deg A_t=1$ in the $G$-grading. We then set
$$
A_t^{\prime}= A_t\vp(e_1)+e_2-\vp(e_1)\mbox{ and }X_t^{\prime\prime}=A_t^{\prime}X_t^{\prime},
$$
where $e_2$ is the identity of $R_2$.

We will first show that
$$
\wD_2=\Sp{X_t^{\prime\prime}\:|\: t\in T}
$$
is a graded subalgebra in $R_2$ isomorphic to $D_1$ (or $D_2$). Now since
$$
\vp(e_1)A_tX_t^{\prime}=\vp(e_1)\vp(X_t)=\vp(X_t)=A_tX_t^{\prime}
$$
and $X_t^{\prime}$ is nondegenerate, it follows that $\vp(e_1)A_t=A_t$. Since $A_t$ and $X_t^{\prime}$ commute in $R_2$, it follows that
$$
X_t^{\prime}A_t\vp(e_1)=\vp(X_t)\vp(e_1)=\vp(X_t)=X_t^{\prime}A_t,
$$
and so
\bee{e996}
\vp(e_1)A_t=A_t\vp(e_1)=A_t.
\ene
In particular, $(e_2-\vp(e_1))A_t=0$.

Now let us recall that $\vp(e_1)=\vp(X_1)=A_1X_1^{\prime}=A_1\in C_{R_2}(D_2)$ and so
\bee{e997}
\vp(e_1)X_t^{\prime}=X_t^{\prime}\vp(e_1)\mbox{ for any }t\in T.
\ene
If we use (\ref{e995}), (\ref{e996}), and (\ref{e997}) we will obtain the following. 
\ben
X_t^{\prime\prime}X_s^{\prime\prime}
 &=& (A_t+e_2-\vp(e_1))X_t^{\prime}(A_s+e_2-\vp(e_1))X_s^{\prime}\\
 &=&  A_tX_t^{\prime}A_sX_s^{\prime}+(e_2-\vp(e_1))X_t^{\prime}X_s^{\prime}\\
 &=& \vp(X_t)\vp(X_s)+(e_2-\vp(e_1))X_t^{\prime}X_s^{\prime}\\
 &=& \vp(X_tX_s)+(e_2-\vp(e_1))X_t^{\prime}X_s^{\prime}\\
 &=& \vp(\alpha(t,s)X_{ts})+(e_2-\vp(e_1))\alpha(t,s)X_{ts}^{\prime}\\
 &=& \alpha(t,s)A_{ts}X_{ts}^{\prime}+\alpha(t,s)(e_2-\vp(e_1))X_{ts}^{\prime}\\
 &=& \alpha(t,s)X_{ts}^{\prime\prime}.
\eqn
Since the elements $X_t^{\prime\prime}$, $t\in T$ are linearly independent, it follows that the mapping $X_t\mapsto X_t^{\prime\prime}$ is a (graded) isomorphism of algebras $D_1$ and $\wD_2$.

Now we denote by $\wC_2$ the centralizer $C_{R_2}(\wD_2)$ of $\wD_2$ in $R_2$. Then $\wC_2$ is a graded subalgebra of $R_2$. The identity element $e_2$ of $R_2$ is in $\wD_2$ since $X_1=e_1$ and
$$
X_{1}^{\prime\prime}=\vp(X_1)+e_2-\vp(e_1)=\vp(e_1)+e_2-\vp(e_1)=e_2.
$$
In this case $R_2=\wC_2\wD_2\cong \wC_2\ot\wD_2$ (see, for instance, \cite{BSZ}).

Now we would like to show that $\wC_2$ is an algebra with elementary $G$-grading. Since $\wC_2$ is a central simple $F$-algebra, by the main result of \cite{BSZ} $\wC_2=\wC_0\wD_0\cong\wC_0\ot\wD_0$ where the grading on $\wC_0$ is elementary while on $\wD_0$ fine. We set $T_0=\su{\wD_0}$. Then $T_0$ is a subgroup in $G$ such that $T_0\cap T=\{ 1\}$. It follows that $\wD_0\wD_2\cong\wD_0\ot\wD_2$ is a graded subalgebra in $R_2$ with fine grading such that $\su{\wD_0\wD_2}=T_0T\cong T_0\times T$. Besides, $R_2\cong\wC_0\ot(\wD_0\ot\wD_2)$ is another decomposition of $R_2$ as the tensor product of algebras with elementary and fine grading. From the proof of \cite[Theorem 6]{BSZ} it follows that the factor $\wD_2$ in the decomposition $R_2=\wC_2\wD_2$ is defined uniquely up to an isomorphism of graded algebras and so $T=\su{\wD_2}=\su{\wD_0\wD_2}=T_0T$. Then $T_0=\{ 1\}$ implying that $\wC_2=\wC_0$ is an algebra with elementary grading.

Now we need to show that $\vp(C_1)\subset\wC_2=C_{R_2}(\wD_2)$. Let $a\in C_1=C_{R_1}(D_1)$. Then $aX_t=X_ta$ for any $t\in T$. Then we have the following.
\ben
\vp(a)X_t^{\prime\prime}
 &=& \vp(a)(\vp(X_t)+e_2-\vp(e_1))\\
 &=& \vp(a)\vp(X_t)+\vp(a)\vp(e_1)(e_2-\vp(e_1))=\vp(aX_t)\\
 &=& \vp(X_ta)=(\vp(X_t)+e_2-\vp(e_1))\vp(a)\\
 &=& X_{t}^{\prime\prime}\vp(a),
\eqn
proving that $\vp(a)\in C_{R_2}(\wD_2)$, that is, $\vp(C_1)\subset\wC_2$.

Finally, let us assume $\vp(e_1)R_2\vp(e_1)=\vp(R_1)$. Since $\vp(C_1)\subset\wC_2$ the containment $\vp(C_1)\subset\vp(e_1)\wC_2\vp(e_1)$ is obvious. To prove the converse, we notice that
$$
\vp(R_1)\cap\wC_2=\vp(e_1)\wC_2\vp(e_1).
$$
Now $C_1=C_{R_1}(D_1)$ and so $b\in R_1$ satisfies $\vp(b)\in C_{R_2}(\wD_2)$ if and only if $b\in C_1$, that is, $\vp(R_1)\cap\wC_2=\vp(C_1)$ which now implies $\vp(C_1)=\vp(e_1)\wC_2\vp(e_1)$.

It remains to look at the homomorphism $\psi:D_1\ra\wD_2$ given by $\psi(X_t)=X_t^{\prime\prime}$. We have the following.
\ben
\vp(e_1)X_t^{\prime\prime}
 &=& \vp(e_1)A_tX_t^{\prime}+\vp(e_1)(e_2-\vp(e_1))X_t^{\prime}\\
 &=& \vp(e_1)A_tX_t^{\prime}=\vp(e_1)\vp(X_t^{\prime})
\eqn
so that $\vp(e_1)\psi(d)=\vp(e_1)\vp(d)$ for any $d\in D_1$ and thus $$\vp(a)\psi(d)=\vp(a)\vp(e_1)\vp(d)=\vp(a)\vp(d).$$ Now the proof is complete.
\qed   

To describe the elementary gradings on infinite-dimensional simple algebras we first consider the case of one \fd matrix algebra embedded in another, both having elementary gradings. Notice that in the claims to follow the grading group $G$ may be infinite and nonabelian.

To start with we notice that if $R=M_n(F)$ is an algebra with an elementary grading given by a tuple $(h_1\ld,h_n)$, $n=km+r$ for some $k,m\ge 1$, $r\ge 0$ and $h_1,\ld,h_n$ satisfy the conditions
\bee{e99e1}
h_{i+1}^{-1}h_{i+2}=h_{i+k+1}^{-1}h_{i+k+2}=\cd=h_{i+(m-1)k+1}^{-1}h_{i+(m-1)k+2}\mbox{ for }1\le i\le k-2
\ene
then the subalgebra $C$ consisting of all block-diagonal matrices of the form $\mathrm{diag}\:\{X$,$\ld$, $X,0\}$ where $X$ is an arbitrary $k\times k$-matrix repeated $m$ times on the diagonal is $G$-graded and isomorphic to a matrix algebra $M_k(F)$ with an elementary grading given by the tuple $(h_1\ld,h_k)$. This easily follows because by (\ref{e99e1}) all matrix units $E_{\alpha+ik,\beta+ik}$, $i=0,1,\ld,m-1$, have the same degree for fixed $1\le\alpha,\beta\le k$. We would like now to prove that any embedding of simple algebras with elementary gradings amounts to this construction.

Let us recall that if $V=\bigoplus_{g\in G}V_g$ is a $G$-graded space then $R=\mathrm{End}\: V$ canonically becomes $G$-graded if, given a $G$-graded basis $\{ v_1,\ld,v_n\}$ of $V$ with $\deg v_i=g_i^{-1}$, $1\le i\le n$, one gives the matrix unit $E_{ij}$ the degree equal $g_i^{-1}g_j$. Thus the grading of $M_n(F)$ induced from $V$ is elementary.

\begin{lemma}\label{l993}
Let $V$ be an $n$-dimensional $G$-graded space over a field $F$ and $\mathrm{End}\: V=R=\bigoplus_{g\in G}R^{(g)}$ the algebra of all linear transformations of $V$ with induced elementary grading. Let $C$ be a graded subalgebra in $R$ which is isomorphic to the matrix algebra $M_k(F)$ with an elementary grading given by the tuple $(g_1,\ld,g_k)$. Then $V$ splits as the sum of $C$-invariant subspaces 
\bee{e99e2}
V=V_1\oplus\cd\oplus V_m\oplus V_0
\ene
where $\dim V_1=\cd=\dim V_m=k$, $V_1,\ld,V_m$ are faithful irreducible $C$-modules while $CV_0=\{ 0\}$. Besides, there is a homogeneous basis of $V$ in which all matrices of the transformations in $C$ have the block-diagonal form $\mathrm{diag}\:\{X,\ld,X,0\}$ where $X$ is a $k\times k$-matrix and the tuple $(h_1,\ld,h_n)$ giving the induced elementary grading on $R=M_n(F)$ satisfies (\ref{e99e1}).
\end{lemma}

\pp Since the grading on $C\cong M_k(F)$ is elementary, any subspace spanned by a set of matrix units is graded. In particular, this is true for any minimal left ideal spanned by all matrix units in a fixed column. Let $L$ be one of such minimal ideals, corresponding to the last, $k$-th column of $C$. If we fix any $v\in V$ then the left $C$-module $Lv$ is either irreducible or equal zero. Moreover, if $v$ is homogeneous then the $C$-submodule $Lv$ is also $G$-graded. These remarks are sufficient to prove the existence of the decomposition (\ref{e99e2}).

Now let $E_{ij}$, $1\le i,j\le k$ be the set of all matrix units of $C$. Since $V_1$ in (\ref{e99e2}) is a faithful $C$-module, there exists a homogeneous element $v\in V$ such that
$E_{1k}v\neq 0$. In this case the vectors $v_i=E_{ik}v=E_{i,i+1}\cd E_{k-1,k}v$, $i=1,\ld,k-1$, $v_k=v$, form a homogeneous basis of $V_1$ and the elementary grading on $C$ induced from this grading is given by the tuple $(g_1,\ld,g_k)$. Indeed, if $\deg v = h$ then $\deg v_i=g_i^{-1}g_kh$, $\deg v_j=g_j^{-1}g_kh$ and so $\deg E_{ij}=g_i^{-1}g_j$ and still $E_{ij}v_j=v_i$. If we choose the bases in other $V_2,\ld,V_m$ and an arbitrary homogeneous basis in $V_0$ then we obtain a realization of $C$ by the the block-diagonal matrices of the form $\mathrm{diag}\:\{X,\ld,X,0\}$.

It remains to consider the tuple $(h_1,\ld,h_n)$ which defines the elementary grading on $R=M_n(F)$ induced from the graded basis of $V$ just constructed. If we denote by $\widetilde{E}_{st}$, $1\le s,t\le n$, the matrix units of $R$ corresponding to this basis, then, as usual, $\deg \widetilde{E}_{st}=g_s^{-1}g_t$. Also, for any $1\le i,j\le k$ we will have
$$
E_{ij}=\widetilde{E}_{ij}+\widetilde{E}_{i+k,j+k}+\cd+\widetilde{E}_{i+(m-1)k,j+(m-1)k}
$$
in $R$ and $\deg E_{ij}=g_i^{-1}g_j$ in $C$, hence in $R$, since the embedding of $C$ in $R$ is graded. Now all $\widetilde{E}_{st}$ are homogeneous and so the conditions (\ref{e99e1}) must be satisfied. Now the proof is complete. \qed

\begin{example} The condition of $C$ having an elementary grading is essential. For example, suppose $C\cong M_n(F)$ with a fine grading. Let $V$ be $C$ itself as a graded vector space and let us assume that $C$ acts on itself by multiplication on the left. Then $R=\ed{V}$ is an algebra with elementary grading induced from $V$ and $C$ a graded subalgebra. So, $C$ is a graded matrix subalgebra of a matrix algebra with an elementary grading but the grading of $C$ is not elementary. So the conclusion of the previous lemma cannot hold for $C$.
\end{example}

Lemma \ref{l993} enables one to describe the gradings on all possible direct limits of matrix algebras with elementary gradings. Here we will need a special case where $C=\bigcup_{i\ge 1}C_i$ where $C_1\subset C_2\subset\cd$ is an ascending chain of matrix algebras and $C_i=e_iC_je_i$, for any $1\le i\le j$, where $e_i$ is the identity element of $C_i$. To start with we generalize the notion of the elementary grading to the case of finitary matrices.

\begin{definition}\label{d991}
Let $R$ be the algebra of finitary matrices and $\emph{\textbf{g}}=(g_1,g_2,\ld)$ a sequence of elements in a group $G$. Then a grading $R=\bigoplus_{g\in G}R^{(g)}$ is called elementary defined by $\emph{\textbf{g}}$ if $R^{(g)}=\Sp{E_{ij}\:|\: g_i^{-1}g_j=g}$.
\end{definition}

\begin{lemma}\label{l994}
Let $C=\bigoplus_{g\in G}C^{(g)}$ be a $G$-graded algebra over a field $F$ which is the union $C=\bigcup_{i\ge 1}C_i$ of an ascending chain of graded matrix subalgebras of orders $n_1$, $n_2$,\ld, with identity elements $e_1,e_2,\ld$. Suppose all the gradings on the subalgebras $C_1, C_2,\ld$ are elementary and $C_i=e_iC_je_i$ for all $i,j$ with $1\le i\le j$. Then $C$ is isomorphic to the algebra $R$ of finitary matrices with elementary grading given by a sequence $\emph{\textbf{g}}=(g_1,g_2,\ld)$ of elements of $G$ in which every $C_i$ is embedded as a graded subalgebra of all matrices with zero entries in all rows and columns whose numbers are greater than $n_i$, $i=1,2,\ld$ The $G$-grading on $C_i$ is elementary given by an $n$-tuple $(g_1,\ld,g_{n_i})$.
\end{lemma}

\pp By Lemma \ref{l993}, we may adjust our graded embeddings in the sequence $C_1\subset C_2\subset\ld$ in such a way that each $C_i$ can be viewed as a graded subalgebra of $C_{i+1}$ consisting of all $n_i\times n_i$ matrices in the left upper corner. These adjustments do not change the isomorphism class of the limit since this depend only on the module structure of $C_{i+1}$ over $C_i$, for each $i$ (see \cite{BZh}). But then the set of all matrices $L_i$ in $C_i$ with zeros outside the first column is a graded subspace of the similar subspace $L_{i+1}$ in $C_{i+1}$. If $\{ e=g_1^{-1}, g_2^{-1},\ld,g_{n_i}^{-1}\}$ is the set of degrees of the matrix units spanning $L_i$ then the elementary grading of $C_i$ is defined by the tuple $(g_1,\ld,g_{n_i})$. Then the set of degrees of the matrix units in $L=\cup_{i=1}^{\infty}$ is the desired sequence of elements of $G$ defining the elementary grading on the algebra of finitary matrices $C$.\qed

\section{Gradings on simple algebras with minimal one sided ideals}\label{sSAFLT}

In this section we consider the gradings by finite abelian groups on simple algebras with minimal one sided ideals. Suppose that $R$ is such an algebra. According to the Structure Theorem in \cite[Chapter 4, Section 9]{NJSR}, there exists a pair of mutually dual spaces $V$ and $\Pi\subset V^\ast$ such that $R\cong V\o \Pi$ with the product is given by 
$$
(v_1\o\pi_1)(v_2\o\pi_2)=\pi_1(v_2)(v_1\o\pi_2)
$$
where $v_1,v_2\in V$, $\pi_1,\pi_2\in \Pi$ and the kernel of the bilinear mapping $(v,\pi)\mapsto \pi(v)$ is trivial. If $\dim V=\dim\Pi=n<\infty$ we have $R\cong M_n(F)$, the matrix algebra of order $n$ over $F$.

The linear mapping $S:V\ra V$ and $T:\Pi\ra\Pi$ are called \emph{adjoint} if $(T(\pi))(v)=\pi(S(v))$. Actually, $T$ is completely defined by $S$ and we write $T=S^\ast$. The Isomorphism Theorem \cite[Chapter 4, Section 11]{NJSR} describes the automorphisms of $V\o \Pi$ with the help of the automorphisms of $V$ in the following way. If $\vp\in \Aut{(V\o\Pi)}$ then there exists a linear automorphism $S: V\ra V$, for which there exists the adjoint automorphism $S^\ast:\Pi\ra\Pi$, such that
$$
\vp(v\o\pi)=S^{-1}(v)\o S^\ast(\pi)\mbox{ for any }v\in V,\,\pi\in\Pi.
$$ 
The automorphism $S$ is defined by $\vp$ uniquely up to a nonzero scalar multiple. 

The finite-dimensional subspaces $V^\prime\subset V$ and $\Pi^\prime\subset \Pi$ are called \emph{compatible} if they are of the same dimension $n$ and the annihilator of $V'$ in $\Pi'$ is zero. As mentioned above, in this case $V^\prime\o\Pi^\prime\cong M_n(F)$. A simple remark is that $V^\prime\o\Pi^\prime\subset  V^{\prime\prime}\o\Pi^{\prime\prime}$ if and only if $V^\prime\subset  V^{\prime\prime}$ and $\Pi^\prime\subset  \Pi^{\prime\prime}$. It is shown in \cite[Chapter 4, Section 16]{NJSR} that $R=V\o\Pi$ has a local system of matrix subalgebras of such form. It will be convenient to label the subalgebras in this local system by the elements of a directed set $I$, that is, an ordered set such that for any $\alpha,\beta\in I$ there is $\gamma\in I$ with $\alpha\prec\gamma$ and $\beta\prec\gamma$. We will have $V_{\alpha}\o\Pi_{\alpha}\subset V_{\beta}\o\Pi_{\beta}$ if and only if $V_{\alpha}\subset V_{\beta}$ and $\Pi_{\alpha}\subset \Pi_{\beta}$. The latter holds if and only if $\alpha\prec\beta$. 

Our aim is to prove the following.

\begin{theorem}\label{p9991} Let a simple algebra $R$ with minimal one sided ideals over an algebraically closed field of characteristic zero be given a grading by a finite abelian group $G$. Then $R$ has a local system of graded finite-dimensional matrix algebras.
\end{theorem}

\pp According to Litoff's Theorem \cite[Chapter 4, Section 15]{NJSR} $R$ is locally matrix, that is, there  a local system $\{ V_{\alpha}\o\Pi_{\alpha}\:|\:\alpha\in I\}$ of matrix subalgebras in a $G$-graded algebra $R=V\o\Pi$. We need to prove that there is another local system   whose terms are $G$-graded matrix subalgebras of the form $\{ \overline{V}_{\alpha}\o\overline{\Pi}_{\alpha}\:|\:\alpha\in I\}$. 

Now the conditions imposed on the field allow one to replace the graded subspaces by the invariant subspaces with respect to the automorphisms corresponding to the multiplicative characters $\chi\in\wg$, given by $\chi\ast r=\chi(g)r$, for any $r\in R^{(g)}$. As mentioned before, to each such $\chi$ one can associate an automorphism $S_\chi:V\ra V$ and its adjoint $S_\chi^\ast:\Pi\ra\Pi$ so that $\chi\ast(v\o\pi)=S_\chi^{-1}(v)\o S_\chi^\ast(\pi)$, for any $v\in V$ and $\pi\in \Pi$. Since $S_\chi$ is defined up to scalar, the mappings $\chi\mapsto S_{\chi}$ and $\chi\mapsto S_{\chi}^{\ast}$ are projective representations of $\wg$ by linear transformations of $V$ and $\Pi$. It is obvious that given a $\wg$-invariant subspace $U$ in $V$, annihilator $U^\perp$ in $\Pi$ is also $\wg$-invariant. 
With these facts in mind, we first pick, for each $\alpha\in I$ a subspace of finite codimension $\Pi_\alpha^\perp$. Set 
$$
U_\alpha=\bigcap_{\chi\in\wg}S_\chi(\Pi_\alpha^{\perp})\subset \Pi_\alpha^\perp.
$$
This is a $\wg$-invariant subspace in $V$ of finite codimension.
Since $\cup_{\alpha\in I} \Pi_\alpha=\Pi$ we must have $\cap_{\alpha\in I} \Pi_\alpha^\perp=0$, hence $\cap_{\alpha\in I} U_\alpha=0$. Note that $U_\gamma
\subset U_\beta$ as soon as $\beta\prec \gamma$. Let us now consider a finite-dimensional $\wg$-invariant subspace $\wg(V_\alpha)=\sum_{\chi\in \wg}\chi(V_\alpha)$. Then there exists $U_\beta$ such that $\wg(V_\alpha)\cap U_\beta=0$. Since $I$ is a directed set, there is $\gamma\in I$ such that $\alpha, \beta\prec\gamma$ hence $\wg(V_\alpha)\cap U_\gamma=0$. Since the projective representation of a finite group is fully reducible there is a $\wg$-invariant subspace $L$ in $V$ such that $V=L\oplus (\wg(V_\alpha)\oplus U_\gamma)=0$. We then set $\overline{V}_\alpha=L\oplus \wg(V_\alpha)$. Also, we set $\overline{\Pi}_\alpha=U_\gamma^\perp$. Since $U_\gamma\subset U_\alpha\subset \Pi_\alpha^\perp$, we have that $\Pi_\alpha\subset\overline{\Pi}_\alpha$. Being an orthogonal complement to a $\wg$-invariant space, $\overline{\Pi}_\alpha$ is $\wg$-invariant. By construction, $\dim \overline{\Pi}_\alpha=\dim \overline{V}_\alpha$ and also $ \overline{\Pi}_\alpha^\perp=U_\gamma$ has trivial intersection with $\overline{V}_\alpha$. This proves that $\overline{V}_\alpha$ and $\overline{\Pi}_\alpha$ are compatible invariant spaces so that $\overline{V}_\alpha\o\overline{\Pi}_\alpha$ is a $\wg$-invariant matrix subalgebra. Since $V\o\Pi_\alpha\subset \overline{V}_\alpha\o\overline{\Pi}_\alpha$, we obtain a $\wg$-invariant local system, hence a local system of graded matrix subalgebras.\bx

\section{Gradings on simple algebras of finitary matrices}

Now we are ready to prove our main result.

\begin{theorem}\label{t991} Let $G$ be a finite abelian group, $R=\bigoplus_{g\in G}R^{(g)}$ be a $G$-graded algebra of infinite matrices each having only finitely many nonzero entries over an algebraically closed field $F$ of characteristic zero. Then $R$ is isomorphic to a graded tensor product $C\otimes D$ where $C$ is such with an elementary grading and $D=M_n(F)$ is a matrix algebra of order $n$ with a fine grading. Additionally, we have $\su{C}\cap\su{D}=\{ 1\}$.
\end{theorem}

\pp According to \cite[Chapter 4, Section 15]{NJSR} $R$ is the same as the simple algebra with minimal one sided ideals in the case where $\dim R$ is countable. Clearly, in this case we can remove unnecessary terms from the local system provided by Theorem \ref{p9991} and conclude that $R$ is the union of the ascending chain $R_1\subset R_2\subset\cd$ of graded simple finite-dimensional subalgebras. Each $R_i$ decomposes as the tensor product $R_i=C_iD_i\cong C_i\ot D_i$ of a simple subalgebra $C_i$ with an elementary grading and a simple subalgebra $D_i$ with a fine grading. The support $T_i=\su{D_i}$ is a subgroup in $G$. Since the number of subgroups in $G$ is finite, by excluding unnecessary subalgebras $R_i$ we may assume that $T_1=T_2=\cd$ is the same subgroup $T$ of $G$. In particular, $\dim D_i=|T|$ for all $i$ and that the $D_i$ as ungraded algebras all isomorphic to the same $M_n(F)$. Since the number of different fine gradings is also finite, we may, as before, assume that all $D_i$ are isomorphic as graded algebras.

Let $\vp_{i+1,i}$ be the graded embedding of $R_{i}$ in $R_{i+1}$. We also set $\vp_{ji} = \vp_{j,j-1}$ $\cd$ $\vp_{i+1,i}$. If we apply Theorem \ref{l992} to each embedding $R_i\subset R_{i+1}$ then we may assume that $\vp_{i+1,i}(C_i)\subset C_{i+1}$ and $\vp_{i+1,i}(C_i)=\vp_{i+1,i}(e_i)C_{i+1}\vp_{i+1,i}(e_i)$ where $e_i$ is the identity element of $C_i$ and that there is an isomorphism $\psi_{i+1,i}(D_i)\ra D_{i+1}$ such that
\bee{e99f1}
\vp_{i+1,i}(a)\vp_{i+1,i}(d)=\vp_{i+1,i}(a)\psi_{i+1,i}(d)\mbox{ for all } a\in C_i,\: d\in D_{i}.
\ene
We set $\psi_1=\mathrm{id}_{D_1}$ and $\psi_i=\psi_{i,i-1}\cd\psi_{2,1}:D_1\ra D_i$, for all $i\ge 2$. Then $\psi_j(d)=\psi_{j,i}(\psi_i(d))$ for any $i,j$ with $1\le i\le j$, and any $d\in D_1$. Besides, using (\ref{e99f1}), we may write
\bee{e99f2}
\vp_{j,i}(a)\vp_{j,i}(d)=\vp_{i,j}(a)\psi_j(\psi_i^{-1}(d))\mbox{ for all } a\in C_i,\: d\in D_{i}.
\ene 
Let us set $C=\bigcup_{i\ge 1}C_i$ and construct an isomorphism  $\rho: R\ra C\ot D_1$. If $a\in C_i$, $d\in D_i$ then we set
\bee{e99f3}
\rho(ad)=a\ot \psi_i^{-1}(d))\mbox{ for any } a\in C_i,\: d\in D_{i}.
\ene
Clearly, (\ref{e99f3}) defines an injective homomorphism of $R_{i}=C_iD_i$ into $C\ot D_1$. Actually, the same formula defines an isomorphism of $R$ to $C\ot D_1$. To prove this we only need to check that $\rho$ is well defined on $R$. Indeed, if  $a\in C_i,\: d\in D_{i}$ and $i<j$ then $\vp_{j,i}(ad) = \vp_{j,i}(a)\vp_{j,i}(d)$ in $R$ and $a=\vp_{j,i}(a)$ in $C$ since we identify $a\in C_i$ with its image $\vp_{j,i}(a)$ in $C_j$. But then, according to (\ref{e99f2}) we should have
\ben
\rho(\vp_{j,i}(a)\vp_{j,i}(d))&=&\rho(\vp_{i,j}(a)\psi_j(\psi_i^{-1}(d)))\\
&=&\vp_{j,i}(a)\ot \psi_i^{-1}(d)=a\ot \psi_i^{-1}(d),
\eqn
proving that, indeed, $\rho$ is defined correctly.

By Lemma \ref{l994} $C$ is isomorphic to the algebra of finitary matrices with an elementary $G$-grading. Since $\su{C}=\bigcup_{i\ge 1}\su{C_i}$ and $T\cap\su{C_i}=\{ 1\}$, for all $i\ge 1$, we have $T\cap\su{C}=\{ 1\}$, and the proof is complete.
\qed

\section{The uniqueness theorem for the elementary gradings of simple algebras of finitary matrices}

The defining sequence $\textbf{g}$ of an elementary grading is not
defined uniquely. In what follows we prove a theorem that gives
necessary and sufficient conditions for two sequences to define
isomorphic gradings. It will be convenient to denote such sequence
as a function $\tau:I\ra G$ such that $\tau(i)=g_i$. Here $I$ is
the sequence of natural numbers or any initial segment of this. In
the latter case we simply deal with $R=M_n(F)$ for a natural
number $n$. With each such function we associate a function
$S_{\tau}:G\ra \mathbb{N}\cup\{\infty\}$ given by
$S_{\tau}(g)=\mathrm{Card}(\tau^{-1}(g))$.

Further notice that for each elementary grading defined by a
function $\tau$ there is a graded vector space $V$ with a basis
$\{ v_i\:|\:i\in I\}$ such that $\deg v_i=g_i^{-1}$. We denote the
subspace spanned by all $v_i$ with $\tau(i)=g$ by $V_{g^{-1}}$. In
this case the algebra of finitary matrices can be identified with
the set all linear transformations of $V$ spanned by the linear
transformations with matrices $E_{ij}$ with respect to the above
basis. The homogeneous component $R^{(g)}$ is then the set of all
linear transformations $\vp$ such that $\vp(V_h)\subset V_{gh}$.

\begin{theorem}\label{tu}
Let $G$ be a group, $R$ and $R^{\prime}$ the algebras of finitary
matrices endowed by two elementary gradings $R=\bigoplus_{g\in
G}R^{(g)}$ and $R^{\prime}=\bigoplus_{g\in G}(R^{\prime})^{(g)}$ defined by
the tuples $\tau$ and $\tau^{\prime}$, respectively. Then $R$ and
$R^{\prime}$ are isomorphic as graded algebras if and only if
there is an element $g_{\,0}\in G$ such that $S_{\tau}(g)=
S_{\tau^{\prime}}(g_{\,0}g)$, for all $g\in G$.
\end{theorem}

\pp First we assume that the gradings defined by $\tau$ and
$\tau^\prime$ are isomorphic. Note that two sequences
$\emph{\textbf{g}}=(g_1,g_2,\ld)$ and
$\emph{\textbf{h}}=(ag_1,ag_2,\ld)$ defines the same gradings on
$R$ and $S_{\tau}(g)=S_{\rho}(ag)$ for all $g\in G$ where
$\rho(i)=ag$. Hence we can suppose that $\rho(1)=e$ that is
$g_1=1$ in $\emph{\textbf{g}}$.

Let $f:R\ra R^{\,\prime}$ be the graded isomorphism of $R$ and
$R^{\,\prime}$, that is, $f(R^{(g)})=(R^{\,\prime})^{(g)}$, for all $g\in G$. Let
us consider the identity components $R^{(e)}$ and $(R^{\:\prime})^{(e)}$.
Each of these algebras is the sum of simple ideals $M^{(g)}$ and
$(M^{\prime})^{(g)}$ each defined as the linear span of the set of
matrix units $E_{ij}$ or $E_{ij}^{\:\prime}$, respectively, such
that $\tau(i)=\tau(j)=\tau^{\prime}(i)=\tau^{\prime}(j)=g$.

Since $f(R^{(e)})=(R^{\:\prime})^{(e)}$ we must have
$f(M^{(g)})=(M^{\prime})^{(\sigma(g))}$ for a bijective map
$\sigma:\su{R}\mapsto \su{R}^{\prime}=\su{R}$ on $G$. Let us also
recall \cite[Corollary 2, Section 4.11]{NJSR} that there is a
linear bijective map $\alpha:V\ra V^{\:\prime}$ such that
$f(\vp)=\alpha\vp\iv{\alpha}$ for any $\vp\in R$. Let us notice
first that such $\alpha$ must satisfy the equation
$\alpha(V_{g^{-1}})=V_{\sigma(g)^{-1}}^{\:\prime}$. Indeed, we
have
$$
M^{(g)}=\{\vp\in R^{(e)}\:|\:\vp(V)\subset V_{g^{-1}}\}\mbox{ and
}(M^{\prime})^{(g)}=\{\vp\in
(R^{\:\prime})^{(e)}\:|\:\vp(V^{\:\prime})\subset
V_{g^{-1}}^{\:\prime}\}.
$$
We have $\alpha M^{(g)}\alpha^{-1}=(M^{\prime})^{(\sigma(g))}$ and
so $\alpha M^{(g))}=(M^{\prime})^{(\sigma(g))}\alpha$. Applying
both sides to $V$ and having in mind the equations
$$
\alpha(V)=V^{\:\prime},\; M^{(g)}(V)=V_{g^{-1}}\mbox{ and }
(M^{\prime})^{(\sigma(g))}(V^{\:\prime})=V_{\sigma(g)^{-1}}^{\:\prime}
$$
we obtain $\alpha(V_{g^{-1}})=V_{\sigma(g)^{-1}}^{\:\prime}$.

Now let us use $\alpha R^{(g)}\alpha^{-1}=(R^{\:\prime})^{(g)}$ or
$\alpha R^{(g)}=(R^{\:\prime})^{(g)}\alpha$, for all $g\in\su{R}\subset
G$. Applying both sides of this equation to any $V_h$,
$h\in\su{V}\subset G$, we obtain $\alpha
R^{(g)}(V_h)=(R^{\:\prime})^{(g)}\alpha(V_h)$ and so
$\alpha(V_{gh})=(R^{\:\prime})^{(g)}(V_{\sigma(h^{-1})^{-1}}^{\:\prime})$.
In other words, $V_{\sigma(h^{-1}g^{-1})^{-1}}^{\:\prime}=
V_{g\sigma(h^{-1})^{-1}}^{\:\prime}$ and \bee{e993c}
{\sigma(h^{-1}g^{-1})^{-1}}=g\sigma(h^{-1})^{-1} \ene for any
$h\in \su{R}, g\in\su{R}$. Recall that $g_1=e$ in
$\emph{\textbf{g}}$ that is $e^{-1}=e\in \su{V}$. Substituting
$h=e$ in (\ref{e993c}) and setting $g_{\,0}=\sigma(e)$, we obtain
$\sigma(g^{-1})=g_{\,0}g^{-1}$, for any $g\in \su{R}$. Note that
for any elementary grading oh $g^{-1}\in \su{R}$ if and only if
$g\in \su{R}$. Hence also $\sigma(g)=g_{\,0}g$ for all $g\in
\su{R}$. So we have $\dim V_g=\dim V_{g_{\,0}g}^{\:\prime}$. Since
$S_{\tau}(g)=\dim V_g$, we easily obtain the desired condition:
there is $g_{\,0}\in G$ such that
$S_{\tau}(g)=S_{\tau^{\prime}}(g_{\,0}g)$ for all $g\in G$.

To prove the converse, we consider two $G$-graded finitary matrix
algebras $R=\oplus_{g\in G}R^{(g)}$, $R^{\:\prime}=\oplus_{g\in
G}(R^{\:\prime})^{(g)}$ and assume that there is $g_{\,0}\in G$ such
that $S_{\tau}(g)=S_{\tau^{\prime}}(g_{\,0}g)$, for any
$g\in\su{R}\subset G$.

We define an isomorphism $f:R\ra R^{\,\prime}$ in the following way.
For each $g\in G$, let the ordered subset $I_g$ label the elements
$v_i$ of the basis of $V\cap V_{g^{-1}}$. Let $I_g^{\prime}$ be
the same thing for $V^{\:\prime}$. Then there is an ordered map
$\beta_g:I_g\ra I_{g_{\,0}g}^{\prime}$. We extend it to a
bijection $\beta$ of $I$ into itself. Then $\beta$ satisfies the
following condition. If $\emph{\textbf{g}}=(g_1,g_2,\ld)$ and
$\emph{\textbf{h}}=(h_1,h_2,\ld)$ then
\begin{equation}\label{e993d}
h_{\beta(i)}= g_{\,0}g_i.
\end{equation}
Denote by $f$ the linear map $R\ra R^{\prime}$ such that
$f(E_{ij})=E_{\beta(i)\beta(j)}$. Then $f$ is an isomorphism and
$f(R^{(g)})= (R^{\:\prime})^{(g)}$ due to (\ref{e993d})
\qed

\begin{remark}\label{r991} The theorem above is no longer true if
we replace the algebra of finitary matrices by the other direct
limits of matrix algebras. For example, suppose an algebra $R$ is
the direct limit of the algebras $R_i=M_{2^i}$, $i=1,2,\ld$, where
the structure mappings $\vp_i:R_i\ra R_{i+1}$ are given by
$X\mapsto\mathrm{diag}\{ X,X\}$. Then the elementary grading by
$G=\langle a\rangle_2$ of $R_i$ given by a tuple $\tau$ can be
extended to the grading of $R_{i+1}$ defined by $\tau^{\prime}$ to
make $\vp_i$ graded if we either choose $\tau^{\prime}=(\tau,\tau)$
or $\tau^{\prime}=(\tau,a\tau)$. If we start wish the grading of
$R_1$ defined by $\tau=(e,a)$ and consider the identity component of
the grading in each of the two cases then we will see the limits of
semisimple algebras, each of which is the sum of two isomorphic matrix
subalgebras. But the Bratteli diagrams \cite{Br} of these limits are
different and so the limits are not isomorphic. At the same time the
``Steinitz numbers'' $S_{\tau}$ and $S_{\tau}^{\:\prime}$ are the same
and both equal to $e^{\infty}a^{\infty}$.
\end{remark}
\begin{remark}\label{r992}
A uniqueness theorem for the $G$-gradings of matrix algebras over
algebraically closed field $F$ of characteristic zero has been
established by A. A. Chasov \cite{Ch}.
\end{remark}

\end{document}